\documentclass[10pt,a4paper]{amsart}
\input amssym.def
\input amssym.tex
\title[ Wohlfahrt's Theorem  ]
{ Wohlfahrt's Theorem and index formula for elementary matrix groups
 and $SL(2, \mathcal O)$ }

\author{ \SMALL  cheng lien  lang}
\author{ \SMALL  mong lung lang}

\begin{document}

\baselineskip=12pt


\thanks {The authors would like to thank K. Petersen for very helpful suggestons
during the preparation of the present article.}
\maketitle

\vspace{-0.3in}

\begin{abstract} The present article determines the indices of
 the principal congruence subgroups of the Bianchi groups $B_d$, $SL(2,\mathcal O )$
 and elementary matrix group $E$  and extends
  Wohlfahrt's Theorem to $B_d$, $SL(2,\mathcal O)$ and $E$,
 where $\mathcal O$ is the ring of integers of some number fields.


\end{abstract}

\section{Introduction}

\subsection{} Elementary matrix groups and $\Gamma = SL(2,\mathcal O)$.
 Let $\mathcal O$ be the ring of integers of a number field
and
   let $\{w_1 = 1, w_2,\cdots w_n\}$ be an integral basis of $\mathcal O$. The elementary matrix group
    studied by P. M. Cohn
    is  the subgroup $E$ of $\Gamma $ generated by (see  [C1], [C2],  and Section 4.2 of [F])
    $$ \left < T_i= T_{w_i}  =
 \left (\begin{array}{cc}
   1 & w_i\\
   0 & 1
   \end{array}\right ),
  S=
 \left (\begin{array}{cc}
   0 & 1\\
   -1 & 0
 \end{array}\right )
 \,:\,  1 \le i \le n
 \right  >.\eqno(1.1)$$
 Let $\pi$ be a nontrivial  ideal of $\mathcal O$.
  The  principal congruence subgroups $E(\pi)$  and $\Gamma (\pi)$ of $E$
   and $\Gamma $ respectively are   defined  by
    $$E(\pi) =\{ (a_{ij})\in E  \,:\,a_{11} - 1, a_{22} -1, a_{12}, a_{21} \in \pi\},\eqno(1.2)$$
    $$\Gamma (\pi) =\{ (a_{ij})\in \Gamma  \,:\,a_{11} - 1, a_{22} -1, a_{12}, a_{21} \in \pi\}.\eqno(1.3)$$
    A subgroup $K$ of $E$  is  called congruence if it contains $E(\pi)$
   for some $\pi \subseteq \mathcal O$.
      Congruence subgroups for $\Gamma $ and $PSL(2, \mathcal O)$ can be defined analogously
     (see  3.3 for $PSL(2, \mathcal O$)).
\subsection{}Index formula.
    The first part of the present article studies the index of the
    congruence subgroup $E(\pi)$ of $E$.
  As a byproduct, the index of the principal congruence subgroup $\Gamma (\pi)$
  of $\Gamma$
    can be recovered $simultaneously$ (see [F], [P] when the class number of $\mathcal O$
     is one).

 \smallskip
    \noindent {\bf Theorem  2.6.}
     {\em Let $(1) \ne \pi \subseteq \mathcal O$ be an ideal and let $ E$ and
      $\Gamma$ be given as in
      subsection $1.1$.
 Then
      $$
       [E : E(\pi)]=[\Gamma   : \Gamma(\pi)]= N(\pi)^3\prod( 1- N(P)^{-2})\eqno(1.4)$$
     }

  Theorem 2.6 is equivalent to $g: E \to SL(2, \mathcal O/\pi)$
 defined by $ g(X) =X$ (mod $\pi)$ is surjective. Unlike the modular group case (Lemma 1.38 of [Sh]),
 we give an indirect proof by studying the  normal series of $ E(\pi)\triangleleft
  E$ that consists of principal congruence subgroups
  (Section 2). Our study shows that $[\Gamma : \Gamma (\pi)]= [E : E(\pi)]$ is
   multiplicative (Lemma 2.1).
    In the case $\pi = \mathcal P^n$, where $\mathcal P$
     is a prime ideal,
     $  \Gamma /\Gamma(\mathcal P)\cong     E/E(\mathcal P)\cong SL(2, \mathcal O/\mathcal P)$ and
     $  \Gamma (\mathcal P^r)/\Gamma (\mathcal P^{r+1}) \cong  E(\mathcal P^r)/E(\mathcal P^{r+1})$ is  elementary abelian of order $N(\mathcal P)^3$ for
      every $r\ge 1$ (Lemmas 2.3 and 2.4). Consequently, our approach determines the order
       and
       group structure of $\Gamma/\Gamma (\pi) $ and $E/E(\pi)$.

\subsection{} Whlfahrt's Theorem for  $E$ and $\Gamma $.
 Let
 $K$ be a subgroup of $\Gamma$ of finite index.
We say $K$ has level $m$ if $m$ is the smallest positive integer
 such that the normal closure of the subgroup generated by
$$\left \{
 \left (\begin{array}{cc}
   1 & mw_i\\
  0 & 1
 \end{array}\right )
 \,:\, 1\le i\le n
 \right \}\eqno(1.5)$$
is contained in $K$. Note that such $m$ exists as $K$ is of finite index in $\Gamma $.
The level of subgroups of $E$ is defined similarly. In [W], Wohlfahrt proved that every subgroup $S$ of finite index of the modular group $\Gamma = SL(2, \Bbb Z)$ is congruence if and only if $\Gamma (m) \subseteq S$ where $m$ is the level of $S$.
The second part of the present article is to prove  the following theorems which
 extend Wohlfahrt's Theorem  to  $E$ and $\Gamma =SL(2, \mathcal O)$.

\noindent  {\bf Theorem 3.2.} {\em Let
     $K$ be a subgroup of  finite index of $E$. Suppose
      that $K$ has level $m$. Then $K$ is congruence if and only of
     $E(m) \subseteq K.$}

\noindent  {\bf Theorem 3.4.} {\em Let
     $K$ be a subgroup of  finite index of $\Gamma$. Suppose
      that $K$ has level $m$. Then $K$ is congruence if and only of
     $\Gamma(m) \subseteq K.$}

\noindent
Theorem 3.4  has been proved by Fine and Petersen  when $\mathcal O $
  is the ring of integers of  $\Bbb Q (\sqrt { d})$
   where  $-d = 1,2,3,7,11,19,43,67,163$
  (Theorem 4.7.3 of
 [F], Theorem 3.1.1 of
 [P]).
As an immediate application of Serre's result [S] and Theorem 3.4, one has

\smallskip
\noindent {\bf Corollary 3.5.} {\em Suppose that $\mathcal O$ has a unit of infinite order
 and that $K$ is normal of  index $m$ in $\Gamma = SL(2, \mathcal O)$. Then $\Gamma (m )
  \subseteq K$.}

\noindent     Theorem 3.4 enables us to construct non-congruence subgroups for $SL(2, \mathcal O$)     and $PSL(2, \mathcal O)$.
     Unlike most of the known construction
 which requires the
     existence of a subgroup $S$ of finite index
    such that $S/[S,S]$ has rank 2 or more (see [F], [GS], [Lu], [N],  and [Z] for examples), our construction
     can be achieved  even if
      $S/[S,S]$ is finite. To be more precise,

 \smallskip
 \noindent {\bf Lemma 5.2.} {\em
   Let $S$ be a  subgroup of $PSL(2, \mathcal O)  $ of index $g$, level $n$. Suppose that
   $\mathcal O \ne \Bbb Z$. Suppose further that
   $S/[S, S]$ has a subgroup $N$ of prime index  $q$ such that
   $q$ is inert or $ q\ge 5$ is split, and
     gcd$\,( q,  |SL(2, \mathcal O/n)|/g)$
     $=1$.
    Then $N$ is a non-congruence subgroup of $PSL(2,\mathcal O)$ of  index $gq$.
    In particular, if the rank of $S/[S,S]$ is one or more, then
   $S$ contains  a  non-congruence subgroup of $PSL(2,\mathcal O)$ of finite index. }

\smallskip
See Definition  5.1 for the   terms {\em split, inert} and {\em ramified}.


 \subsection{} Congruence subgroup problem.
  We say $ G =PSL(2, \mathcal O)$
   has congruence subgroup property (CSP) if every subgroup of finite index
    is congruence.
   Applying the well known results of Serre [S],
 $PSL(2, \mathcal O)$ does not have CSP
 if and only if $\mathcal O$ is either $ \Bbb Z$ or  the ring of integers of the
  imaginary quadratic field $\Bbb Q (\sqrt d)$. In the latter case,
    the group is  known as  the Bianchi group $B_d$.
     Equivalently, $ B_d = PSL(2, O_d)$, where $O_d$ is
     the ring of integers of the
  imaginary quadratic field $\Bbb Q (\sqrt d)$.
  The third part of the present article is to study the commutator
       and power subgroups of $B_d$.
     By Theorem 3.4, we have

\smallskip
\noindent {\bf Proposition  6.2 and Proposition 6.3.} {\em
Suppose that
$ d \not \equiv 5\,\,(mod\, 8)$. Then
 $ B_d ^2 =\left < x^2 \,:\, x\in B_d\right >$ is non-congruence if and only if $|B_d/B_d^2| \ge 8$.
  In particular, if the class number of $O_d$ is three or more, then
   $B_d^2$ is non-congruence. In the case $d\equiv 5\,\,(mod\, 8)$, $B_d$ is non-congruence
    if and only if $d \ne -3$.
  }

     \smallskip
     \noindent {\bf Proposition 5.3.} {\em Let $q\in \Bbb N$ be
a prime.
Suppose that  $q$ is inert or $q\ge 5$ is split in $O_d$  and that $d \ne -1,-3$.
      Then $B_d$
       has a normal non-congruence subgroup of level $q$, index $q$.}

\smallskip Normal non-congruence subgroups for $B_{-1}$ and $B_{-3}$
 can be found by applying Lemma 5.2
 (see Proposition 5.4). Consequently, the fact that $B_d$ does not have
  CSP  can be viewed as a  consequence
   of Wohlfahrt's Theorem.
    Suppose that $G = PSL(2, \mathcal O)$ has CSP.
   It is well known that $S/[S,S]$ is finite if $[G :S]$ is finite (Section 16 of [BMS]). With the help
     of Wohlfahrt's Theorem, one has the following slightly better result.

 \smallskip
 \noindent {\bf Proposition 7.1.} {\em Suppose that $G= PSL(2, \mathcal O)$ has CSP.
     Let $S$ be a  subgroup of $G$ of index $g$, level $n$.
      Then $S/[S, S]$ is finite.
       Let $p$ be a  prime divisor  of
      $|S/[S, S]|$.
      Suppose that $p$ is inert or $p\ge 5$ is split. Then $p$
       divides
   $   |SL(2, \mathcal O/n)|/g$.
      In particular, if $p$ is a prime divisor of   $|G/[G,G]| <\infty$,  then $p$ is ramified or split.
 In the case $p$ is split,  $p = 2, 3$. }

 Applying Proposition 7.1, we show that if $G= PSL(2, O_d)$ has CSP, then
  $G/[G,G]$ is of order $2^a3^b$ for some $a, b$ (Proposition 7.2). This makes
   groups with CSP very different from the Bianchi groups $B_d$ as
   $B_d$ do not have CSP and $B_d/[B_d, B_d]$ is finite only of $ d = -1, -3$
    (Theorem 4.1 and Table 1).

\subsection{} Discussion. Our main results (Theorem 2.6 and 3.2. 3.4) are in line with
 the modular group case and their proofs are elementary
  (compare to the works we listed in the references). However,  it seems worth recording
 as both the index formula and Wohlfahrt's Theorem are of
  importance in the study of the Bianchi groups,  $SL(2, \mathcal O)$ and the
   elementary matrix groups.


   \section{Index Formula of   $[\Gamma  : \Gamma (\pi)]$ and    $[E : E(\pi)]$   }
  The  purpose of this section is to determine the indices
   $[\Gamma  : \Gamma (\pi)]$
    and   $[E : E(\pi)]$ {\em simultaneously}.
    Recall that $g : E \to SL(2, \mathcal O/A)$ is defined by
    $g(X) = X$ (mod $A$) and that a  key fact of  the study of the index $[E : E(A)]$
     is to prove that $g$ is surjective (see subsection 1.2). We find this difficult
      as a matrix  $\sigma\in SL(2, \mathcal O)$ satisfies  $\sigma \equiv x$, where $x \in SL(2,
      \mathcal O/A)$,  is not necessarily a member of $E$ as $E \ne SL(2, \mathcal O)$
      if the class number of $\mathcal O$ is two or more.
       As a consequence, we have to prove the
        surjectivity of $g$
        in a different way.

  \subsection{} Upper and  lower bounds for   $[\Gamma  : \Gamma (\pi)]$.  It is clear
     that $$[E : E(\pi)]=
     [\Gamma \cap E : \Gamma (\pi)\cap E]\le[\Gamma  : \Gamma(\pi)].\eqno(2.1)$$
      Hence $[E : E(\pi)]$ is an lower bound of  $[\Gamma  : \Gamma (\pi)]$.
       To start off our study of an upper bound of   $[\Gamma  : \Gamma (\pi)]$, we consider the homomorphism
    $ f\,:\, \Gamma  \to SL(2, \mathcal O/\pi)$ defined by $f(X) = X$ modulo $\pi$.
    It is clear that the kernel of $f$ is
     $\Gamma (\pi)$. As a consequence,
     $$[\Gamma  : \Gamma (\pi)] \le |SL(2, \mathcal O/\pi)|.\eqno(2.2)$$
     Consider the prime factorisation $\pi = \prod \pi_i^{e_i}$ and the homomorphism
      $g \,:\, SL(2, \mathcal O /\pi) \to \prod SL(2, \mathcal O /\pi_i^{e_i})$
      defined by
      $g(X) = (X \mbox{ (mod } \pi_1^{e_1}),\, X \mbox{ (mod } \pi_2^{e_2}), \cdots
      ,   X \mbox{ (mod } \pi_r^{e_r})   ).$
       It is clear
       that $g$ is injective. As a consequence,
   $$[\Gamma  : \Gamma (\pi)] \le |SL(2, \mathcal O /\pi)|
   \le \prod |SL(2, \mathcal O/\pi_i^{e_i})|
   .\eqno(2.3)$$
 Similar to the modular group
    case (see Section 1.6 of [Sh]),
    the order of $ SL(2, \mathcal O /\pi_i^{e_i})$ can be determined.
    It is
  $$ |SL(2, \mathcal O /\pi_i^{e_i})|=
N(\pi_i)^{3e_i}(1-N(\pi_i)^{-2}),\eqno(2.4)$$
\noindent where $ N(\pi_i)$ is the absolute norm of $\pi_i$.
  As a consequence,
  $[\Gamma  : \Gamma (\pi)] \le
     N(\pi)^3 \prod (1-N(\pi_i)^{-2})$. In summary, one has
  $$[E : E(\pi)]\le [\Gamma  : \Gamma (\pi)] \le
     N(\pi)^3 \prod (1-N(\pi_i)^{-2}).\eqno(2.5)$$
   The rest of this section is to study    $[E : E(\pi)]$.
    See Theorem 2.6 for the main results.


\subsection{} The index  $[ E : E(\pi)] $ is multiplicative.
We  prove the following useful lemma which will be used later in our study
 of the concept {\em level} of Wohlfahrt's (Lemma 3.1).

 \smallskip
 \noindent {\bf Lemma 2.1.} {\em Let $K$ be a normal subgroup of $E$.
  Suppose that $T_i\in K$ for all $1\le i \le n$. Then $K = E$.
  In particular, let $A$ and $B$ be two ideals of $\mathcal O $ such  that $A+B = \mathcal O$.
   Then $E(A)E(B)= E$,
     $[E : E(AB) ]  = [E(A):E(AB)][E(B):E(AB)]=[E :E(A)][E:E(B)]$.
   }

  \smallskip
  \noindent {\em Proof.}
   Since $S^4 =1$, $T_1 \in K \triangleleft E$, $(ST_1)^4 \in K$.
         Note that the order of $T_1S$ is three. Hence
          $$T_1S = (T_1S)^4 \in K.\eqno(2.6)$$
           Hence $S, T_i \in K$.  It follows that $K = E$.
            This completes the proof of the first part of the lemma.
          Let $A$ and $B$ be two ideals of $\mathcal O$ and let
  $\{a_1, a_2 \cdots, a_n\}$ and $\{  b_1, b_2, \cdots , b_n\}$ be integral bases of $A$ and $B$
   respectively. Since $A+B = \mathcal O $, there exists $x_{i_k}, y_{j_k} \in \Bbb Z$ such
    that
    $$\sum x_{i_k} a_k + \sum y_{i_k} b_k = w_i\eqno (2.7)$$
     for all $1\le i \le n$.
 Pass to $E(A)$ and $E(B)$, (2.7) implies that (see Remark 2.2 for quadratic case)
       $$T_{w_i} =  T_i  \in E(A) E(B)\mbox{ for all } i. \eqno(2.8)$$
      Apply the first part of our lemma,   $E= E(A)E(B)$.
            Since $E(A)\cap E(B) = E(AB)$, by second isomorphism theorem, we have
            $[E : E(AB) ]  = [E(A):E(AB)][ E(B):E(AB)]
            . $\qed

\smallskip
\noindent {\bf Remark 2.2.}
        Let $A$ and $B$ be two ideals of $\mathcal O$,
        where $\mathcal O = (1, w) $ is the ring of integers of $\Bbb Q(\sqrt d)$ for some $d\in \Bbb Z$
         (square free).
         Let
  $\{a_0, a_1+a_2 w\}$ and $\{ b_0, b_1+ b_2 w\}$ be integral bases of $A$ and $B$
   respectively. Since $A+B = \mathcal O$, there exists $x_i, y_j \in \Bbb Z$ such
    that
    $ x_0a_0 + x_1(a_1+ a_2w) +  x_2 b_0 + x_3( b_1 +b_2w) =1$ and
      $y_0a_0 + y_1(a_1+ a_2w) + y_2b_0 + y_3(b_1 +b_2w)=w$. As a consequence,
      $$T_1^{x_0a_0} (T_1^{a_1}T_w^{a_2})^{x_1}
      T_1^{x_2b_0} (T_1^{b_1}T_w^{b_2})^{x_3} =  T_1,\,
    T_1^{y_0a_0} (T_1^{a_1}T_w^{a_2})^{y_1}
      T_1^{y_2b_0} (T_1^{b_1}T_w^{b_2})^{y_3} =  T_w
    .\eqno(2.9)$$
     Note that $T_1^{a_0}, $ $ T_1^{a_1} T_w^{a_2} \in E(A)$ and that $T_1^{b_0}, $
     $T_1^{b_1}T_w^{b_2} \in E(B)$.
      As a consequence,
       (2.10) implies that
      $ T_1, T_w \in E(A) E(B). $

\subsection{} The index $[ E : E(\pi^m)]$, where $\pi$ is a prime ideal.
 We shall first study the index $[E : E(\pi)]$. It is clear that
  $E/E(\pi)$ is a subgroup of $SL(2, \mathcal O /\pi)$. Note that $\mathcal O/\pi$ is a finite
   field of characteristic $p$ where $p$ is the smallest positive rational prime  in $\pi$.

\smallskip
\noindent {\bf Lemma 2.3.} {\em
Let $\pi\subseteq \mathcal O $ be a prime ideal and let $p$ be the smallest positive
rational   prime in $\pi$.
 Then
$E/E(\pi) \cong SL(2, \mathcal O/\pi)$,
where $\mathcal O/\pi$ is a finite field of characteristic $p$   and  $|SL(2, \mathcal O /\pi)|= N(\pi)^3(1-N(\pi)^{-2})$ .
}

 \smallskip
 \noindent {\em Proof.}
Since $\{w_1, w_2, \cdots , w_n \}$ is an integral basis of $\mathcal O $
(see subsection 1.1) and $S, T_i
\in E$ for all $i$. $E/E(\pi)$ contains all the elementary matrices.
$$S,\,T_x=\left (\begin{array}{cc}
1& x\\
0&1\\
\end{array}\right),\,L_y=
\left (\begin{array}{cc}
1&0\\
y&1\\
\end{array}\right ) \in E/E(\pi)\eqno(2.10)$$
for all $x,y \in \mathcal O/\pi$. Note that
$\mathcal O/\pi$ is a finite field.
 One may  perform elementary row and
column operations to show that every matrix in $SL(2, \mathcal O/\pi)$
can be written as a word in $S, T_x$ and $ Ly$. As a consequence,
$E/E(\pi) \cong SL(2, \mathcal O/\pi)$.
This completes the proof of the lemma.\qed




\smallskip
\noindent {\bf Lemma 2.4.} {\em
Let $\pi\subseteq \mathcal O$ be a prime ideal, $m\ge 1$.
Then $E(\pi^m )/E(\pi^{m+1})$ is an elementary abelian $p$-group of
 order  $   N(\pi)^3$, where $p$ is the smallest positive rational prime in $\pi$.
 }

 \smallskip
 \noindent {\em Proof.} Suppose that  $ N(\pi) = |\mathcal O/\pi| = p^r$ for some $r$.
  It follows that
    $\pi^n/\pi^{n+1}$ is an elementary abelian  $p$-group of order $p^r$.
     Let $\{x_1, x_2, \cdots , x_r\}$ be a set of generators of $\pi^n /\pi^{n+1}$.
     Then $x_i \in \pi^n \setminus \pi^{n+1}$ and
  $ \{ x_1,  x_2, \cdots ,  x_r\}$ is
  a $\Bbb Z_p$-independent set. Equivalently,
  if $\sum_{i=1}^r a_i x_i \in \pi^{m+1}$,
 where $a_i \in \Bbb Z$, then
 $p|a_i$ for all $i$.
  Since $T_i \in E$ for all $i$,
   it is clear
   that
   $$ X_i =
  \left (\begin{array}{ll}
  1 & x_i \\
  0 & 1
\end{array}\right )\in E(\pi^m)\setminus  E(\pi^{m+1})\eqno(2.11)$$
 for all $x_i$. Set
     $ Y_i =  SX_i S^{-1}, $ $ Z_i = T_1S X_iS^{-1}T_1^{-1}$
      (see (1.1) for the definition of $T_i$). The matrix
       forms of $X_i, Y_j, Z_k$ are given as follows.
  $$
  \left (\begin{array}{ll}
  1 & x_i \\
  0 & 1
\end{array}\right ),
\left (\begin{array}{cc}
  1 & 0 \\
 -x_i  & 1
\end{array}\right ),
 \left (\begin{array}{cc}
  1-x_i  & x_i\\
- x_i   & 1+x_i
\end{array}\right )\in E(\pi^m)\setminus  E(\pi^{m+1})\eqno(2.12)
$$
  Write the above matrices into the form $ I + A_i$ for some
   $A_i$. Since $x_ix_j \in \pi^{m+1} $ for all $i$ and $j$,   one sees easily that
  \begin{enumerate}
  \item[(i)]
 $ (I +  A_i)^p \equiv I$ (mod  $\pi^{m+1})$,
 \item[(ii)] $(I+A_i)(I+ A_j )\equiv
 I+(A_i+A_j)$ (mod  $\pi^{m+1}).$
\end{enumerate}
As a consequence, the above  matrices modulo $\pi^{m+1}$ generate an elementary abelian $p$-group. Set
   $$ A = \left <
  \left (\begin{array}{ll}
  1 & x_i \\
  0 & 1
\end{array}\right ),
\left (\begin{array}{cc}
  1 & 0 \\
 -x_i  & 1
\end{array}\right )\right >\,, B= \left <
 \left (\begin{array}{cc}
  1-x_i & x_i\\
- x_i   & 1+x_i
\end{array}\right )\right > \eqno(2.13)$$
Since $ \{x_1,  x_2, \cdots ,x_r\}$ is a
   $\Bbb Z_p$-independent set,
  one sees easily that $A$ modulo $\pi^{m+1}$ is elementary abelian of order $p^{2r}$,
 $B$ modulo $\pi^{m+1}$ is elementary abelian of order $p^{r}$ and
  $A\cap B = \{1\}$  modulo $\pi^{m+1}$ (see  Remark 2.5).
  Hence  $|E(\pi^m /E(\pi^{m+1}) | \ge p^{3r}= N(\pi)^3$.
  Applying (2.5), $|E(\pi^m /E(\pi^{m+1}) | \le N(\pi)^3$.
    This completes the proof of  the lemma.
     \qed

\smallskip
\noindent {\bf Remark 2.5.}
 Since $A$ and $B$ are abelian groups, elements in $A$ and $B$ take the following simple forms respectively.
$$ I +   \sum_{i=1}^r  \left (\begin{array}{cc}
  0 &  a_i  x_i \\
 -   b_i  x_i     & 0
\end{array}\right ),\,\,
 I +
 \sum_{i=1}^r
 \left (\begin{array}{cc}
  - c_i  x_i     & c_i x_i \\
 -    c_i  x_i     &  c_i  x_i
\end{array}\right ).\eqno(2.14)$$

\subsection{} The main result. Apply  (2.5) and Lemmas 2.3 and 2.4,
 we  have the following  results.  It gives the indices in terms of the norms  of ideals.

     \noindent {\bf Theorem 2.6.} {\em Let $1\ne \pi \subseteq \mathcal O$ be
      an ideal. Then
      $$ [\Gamma  : \Gamma (\pi)]=[E : E(\pi)]= N(\pi)^3\prod( 1- N(P)^{-2}),\eqno(2.15)$$
     where the product is over the set of all prime ideals $P$ that divide $\pi$
     and $N(P)$ is the absolute norm of $P$.}

\smallskip
\noindent {\em Proof.}
 Since the index $[E: E(\pi)] $ is multiplicative (Lemma 2.1), one may apply Lemmas 2.3 and 2.4 to conclude
  that
 $[\Gamma  : \Gamma (\pi)]= [E : E(\pi)]= N(\pi)^3\prod( 1- N(P)^{-2})$.  \qed

\smallskip
\noindent {\bf Corollary 2.7.} {\em The homomorphisms
$f\,:\, \Gamma \to SL(2 ,\mathcal  O /\pi)$, $g\,:\, E  \to SL(2, \mathcal O/\pi)$ defined by
$f(X) =X$, $g(X) = X$ modulo $\pi$ are surjective.}

 \section{Wohlfahrt's Theorem }


\subsection{}  Wohlfahrt's Theorem for $E$.
   Denoted by $N( E(mn),  m)$ the normal closure
    of $E(mn)$ and $T_i^m $ $ (1\le i \le n$)  in $E$. The following lemma (Lemma 3.1) shows that
  $ N(E(mn), m )$
  $ = E(m)$. It follows immediately from this lemma that
  Wohlfahrt's Theorem can be extended to the elementary matrix groups (Theorem 3.2).
 \smallskip

  \noindent {\bf Lemma 3.1.} {\em
   $ N(E(mn), m ) = E(m)$.
  }

  \smallskip
  \noindent {\em Proof.}
  Let $x$ be the smallest positive rational integer such that $
   E(x) \subseteq  N( E(mn),  m)  \subseteq E(m).$
  Suppose that $ x> m$.
  Since $m$ is a divisor of $x$, one has
   $ x= mm_0 p$, where $m_0, p \in \Bbb N$ and $p$ is a prime.
Set $q = mm_0$. Then $x = pq$.  Since $E(pq) =  E(x) \subseteq  N( E(mn),  m) $, it follows that
 $$ N( E(pq), q)
  \subseteq N( E(mn), m).\eqno(3.1)$$

\noindent {\bf Case 1.}  gcd$\,(p, q)=1$.  Let $X= N(E(pq), q)$.
 It is clear that $X \subseteq E(q)$.
We consider the group $X E(p)$. Since gcd$\,(p,q)=1$,
$T_i^p\in E(p)$ and $T_i^q\in X$ for all $i$, we have
$T_i\in X E(p)$ for all $i$.
 This implies that $X E(p) = E$ (Lemma 2.1).
 Hence
 $$
  {\small X/E(pq) \times E(p)/E(pq)
   \cong E/E(pq) \cong E(p)/E(pq)\times E(q)/E(pq).}\eqno(3.2)$$
   Note that $X \subseteq E(q)$.
  It is now clear that (3.2) implies that $X =  E(q)$. By (3.1),
   $E(q)=N(E(pq),q) \subseteq  N( E(mn), m)$.
   This contradicts the minimality of $x$. Hence $x=m$.
    Equivalently,  $ N(E(mn), m ) = E(m)$.

    \smallskip
    \noindent {\bf Case 2.} gcd$\,(p,q)\ne1$. It follows that $p|q$.
     By Theorem 2.6, $[E(q) : E(pq)] = N(p)^3$.
      Since $N(E(pq), q) \subseteq E(q)$,
      $$ [ N(E(pq), q): E(pq)]\le N(p)^3.\eqno(3.3)$$
    It is clear that
   $T_i^{q}, $  $  ST_i^{q}S^{-1}, $ $ T_1ST_i^{q}S^{-1}T_1^{-1}\in E(q)$.
  The matrix form of the above  matrices are given as follows.
  $$
  \left (\begin{array}{ll}
  1 & qw_i \\
  0 & 1
\end{array}\right ),
\left (\begin{array}{cc}
  1 & 0 \\
 -qw_i & 1
\end{array}\right ),
 \left (\begin{array}{cc}
  1-qw_i & qw_i\\
-  qw_i  & 1+qw_i
\end{array}\right ).\eqno(3.4)
  $$
 Similar to Lemma 2.4, the above  matrices modulo $E(pq)$ generate
  an elementary abelian $p$-group of order $N(p)^3$. Since these matrices
   are conjugates of $T_i^q$, it follows
    that $[N(E(pq), q)$
    $ :E(pq)] $ is  at least $N(p)^3$.
     By (3.3),
  $[N(E(pq), q) :E(pq)]= N(p)^3 $. Hence $N(E(pq), q) = E(q)$.
By (3.1),
   $E(q) \subseteq  N( E(mn), m)$.
   This contradicts the minimality of $x$. Hence $x=m$.
    Equivalently,  $ N(E(mn), m ) = E(m)$.
This completes the proof of the lemma. \qed
 \smallskip

 \noindent {\bf Theorem 3.2.} (Wohlfahrt's Theorem)
 {\em Let $K$ be a subgroup of $E$ of finite index. Suppose that
  the level of $K$ is $m$. Then $K$ is congruence if and only if
   $E(m) \subseteq K$.}
   \smallskip

   \noindent {\em Proof.} Suppose that $K$ is congruence. Then
    $E(n) \subseteq K$ for some $n$. Hence $E(mn) \subseteq K$.
    Since the level of $K$ is $m$, it follows
     that $N(E(mn), m) \subseteq K$.
      By Lemma 3.1,  $E(m) \subseteq K$. The converse is clear. \qed

\subsection{} Wohlfahrt's Theorem for $\Gamma = SL(2, \mathcal O) $.
 Denoted by $N( \Gamma (mn),   m)$ the normal closure
    of $\Gamma (mn)$ and $ T_i^m$ $(1\le i \le n)$ in $\Gamma$.
The following lemma is essential in our study of Wohlfahrt's Theorem for
$\Gamma $.

 \smallskip

  \noindent {\bf Lemma 3.3.} {\em  $ N(\Gamma (mn), m ) = \Gamma (m)$.
  }

  \smallskip
\noindent {\em Proof.}
Let $X= N(\Gamma (mn), m )$.  An easy study of the group diagram of the following
 six groups
  $ \Gamma (mn) \subseteq X \subseteq \Gamma (m) \mbox { and }
   \Gamma (mn)\cap E(m)\subseteq  X\cap E(m)\subseteq \Gamma (m)\cap E(m)  $
  implies that
\begin{enumerate}
\item[(i)]
$[\Gamma (m) : X] \ge [\Gamma (m) \cap E(m) : X\cap E(m)] = [E(m) : X \cap E(m)]$,
\item[(ii)]
$[X : \Gamma (mn) ] \ge [X \cap E(m) : \Gamma (mn)  \cap E(m)] = [X\cap E(m) : E(mn)]$,
\item[(iii)] $  [\Gamma (m) : \Gamma (mn)] =  [E(m) : E(mn)] $ (Theorem 2.6).
\end{enumerate}
Applying (iii) of the above to (i) and (ii), it follows that the inequalities are
 actually equalities. Hence (i) of the above becomes
$$[\Gamma (m) : X]=  [\Gamma (m) \cap E(m) : X\cap E(m)] = [E(m) : X \cap E(m)].\eqno(3.5)$$
 It is clear that
$N(E(mn), m ) \subseteq X\cap E(m)$. By Lemma 3.1,
$ N(E(mn), m )= E(m)$. Hence the third term of (3.5) is 1.
This implies that
  $ [\Gamma (m) : X]=1$.
   Equivalently, $ \Gamma (m) = N(\Gamma(mn), m )$.  \qed
\smallskip

 Similar to Theorem 3.2, we may extend Wohlfahrt's Theorem to $\Gamma $
as follows.

 \smallskip
 \noindent {\bf Theorem 3.4.} (Wohlfahrt's Theorem)
 {\em Let $K$ be a subgroup of $\Gamma $ of finite index. Suppose that
  the level of $K$ is $m$. Then $K$ is congruence if and only if
   $\Gamma(m) \subseteq K$.}

\smallskip
\noindent {\bf Corollary 3.5.} {\em Suppose that $\mathcal O$ has a unit of infinite order
 and that $K$ is normal of  index $m$ in $\Gamma = SL(2, \mathcal O)$. Then $\Gamma (m )
  \subseteq K$.}

  \smallskip
  \noindent {\em Proof.} Since $\mathcal O$ has a unit of infinite order, every subgroup of $\Gamma$  of finite index is congruence (Serre [S]).
  Since $K$ is normal of index $m$, the level of $K$ is a divisor of $m$.
   It follows from Theorem 3.4 that
    $\Gamma (m ) \subseteq K$. \qed

\smallskip
\noindent {\bf 3.3.} Discussion. Let $G = PSL(2, \mathcal
O)$ or $PE = E/\left <\pm I
\right >$.
 Define the principal congruence subgroup $G(\pi)$ to be
 $G(\pi) = \{ x \in G
 \,:
 \, x \equiv \pm I \,(mod\, \pi)\}$. A subgroup $K$ of $G$ is a congruence subgroup
  if $K$ contains $G(\pi)$ for some $\pi \subseteq \mathcal O$. Let $K$ be a subgroup of $G$
   of finite index. The level of $K$ is the smallest positive integer $m$ such that
   the normal closure of the subgroup generated by
$$\left \{
 \pm \left (\begin{array}{cc}
   1 & mw_i\\
   0 & 1
 \end{array}\right )
 \,:\, 1\le i\le n
 \right \}\eqno(3.6)$$
is contained in $K$.
Theorems 3.2 and 3.4 can be extended easily to $G$ as follows.

   \smallskip
 \noindent {\bf Theorem 3.6.} (Wohlfahrt's Theorem)
 {\em Let $K$ be a subgroup of $G = PSL(2, \mathcal O) $ or $PE$ of finite index. Suppose that
  the level of $K$ is $m$. Then $K$ is congruence if and only if
   $G(m)  \subseteq K$.}

\section{Known results of the Bianchi groups }

 \subsection{}In the case $\mathcal O =O_d$ is the ring of integers of the imaginary quadratic field
   $\Bbb Q(\sqrt d)$, $d <0$.
    $B_d = PSL(2, O_d)$ is known as the Bianchi group.

           \subsection{} Known results about projective elementary matrix groups $PE_d = E_d/\left < \pm I\right >$ and  Bianchi groups $B_d$.
 Let  $\{1, w\}$ be an integral basis of $O_d$. Then $PE_d$ is generated by
$$  T_1 =\left (
\begin{array}{rc}
1 & 1 \\
0 & 1 \\
\end{array}
\right ),
T_w =\left (
\begin{array}{rc}
1 & w \\
0 & 1 \\
\end{array}
\right ),
S= \left (
\begin{array}{rc}
0 & 1 \\
-1 & 0 \\
\end{array}
\right ).\eqno(4.1)$$
In the case $-d = 1,2,3,7,11$, Cohn's results [C1] implies  that $PE_d = B_d$. The group structure for such $B_d$ is completely known (Theorem 4.3.1 of Fine [F]). In the
 case $-d\ne 1,2,3,7,11$, $PE_d = \{a, t, u\,:\, a^2 = (at)^3 = [t, u] =1\}$ (Theorem 4.8.1
  of [F]).  Set $\Omega=\{-1,-2,-3,-7,-11\}$.
  As an easy application of their results, one has the following.

\medskip

\begin{center}
\mbox {Table 1.  $PE_d/[PE_d , PE_d]$ and $PE_d/PE_d^2$}
\end{center}

 $$
 \begin{array}{|r|r|r|r|}
\hline
 \phantom{\Big |}  d & PE_d/[PE_d,PE_d]& PE_d/PE_d^2 & \mbox{Remark}\\
\hline
 \phantom{\Big |}   -1 & \Bbb Z_2 \times \Bbb Z_2  & \Bbb Z_2 \times \Bbb Z_2 &B_d=PE_d   \\
\hline
 \phantom{\Big |}   -2 & \Bbb Z \times \Bbb Z_6 & \Bbb Z_2 \times \Bbb Z_2 &B_d=PE_d    \\
\hline
 \phantom{\Big |}   -3 & \Bbb Z_3  & 1  &B_d=PE_d   \\
\hline
 \phantom{\Big |}   -7 & \Bbb Z \times \Bbb Z_2  & \Bbb Z_2 \times \Bbb Z_2 &B_d=PE_d    \\
\hline
 \phantom{\Big |}   -11 & \Bbb Z \times \Bbb Z_3  &  \Bbb Z_2  &B_d=PE_d   \\
\hline
 \phantom{\Big |}   d\notin \Omega & \Bbb Z \times \Bbb Z_6  &
  \Bbb Z_2 \times \Bbb Z_2  &B_d\ne PE_d   \\
\hline
\end{array}
$$

\medskip
\noindent
where $B_d^2 = \left < x^2 \,:\, x\in B_d \right > $ and
 $PE_d^2 = \left < x^2 \,:\, x\in PE_d \right > $.
We now turn our attention to the case  $ d\notin \Omega=\{-1,-2,-3,-7,-11\}$.
 Works of Swan and others (see Theorem 6.2.2 and  pp. 195 of Fine [F]) imply that

\smallskip
\noindent {\bf Theorem 4.1.} {\em Suppose that $ -d \ne 1, 2, 3 ,7 ,11$.
Then the rank  $r$ of $B_d/[B_d, B_d]$ is finite. Further,
 $ r  \ge {h_d}$, where
  $h_d$ is  the class number of $O_d$.
  In particular,  $\infty > |B_d/B_d^2| \ge  2^r $.}

\section{First application : Non-congruence subgroups of small indices}
The main purpose of this section is to construct normal non-congruence
 subgroups of small indices and levels of the Bianchi groups.
 Throughout the section, $d <0$
   is square free and $B_d(x)$ is the principal congruence subgroup associate to $x$.

 \smallskip
 \noindent {\bf Definition 5.1.}
Consider the prime decomposition $(q) = \mathcal P_1^{e_1}  \mathcal P_2^{e_2}
 \cdots \mathcal P_s^{e_s}  $ in $\mathcal O$.
  $q$ is {\em inert} if $s =1$, $e_1=1$,
 $q$ is {\em split} if $s\ge 2$,  $e_1 = e_2 = \cdots e_s=1$,
 $q$ is {\em ramified} if $e_i \ge 2$ for some $i$.

 \smallskip
 \noindent {\bf Lemma 5.2.} {\em
   Let $S$ be a  subgroup of $PSL(2, \mathcal O)  $ of index $g$, level $n$.
   Suppose that $\mathcal O \ne \Bbb Z$.
   Suppose further that
   $S/[S, S]$ has a subgroup of prime index  $q$
such that   $q$ is inert or $ q\ge 5$ is split, and
     gcd$\,( q,  |SL(2, \mathcal O/n)|/g)$
     $=1$.
    Then $N$ is a non-congruence subgroup of $PSL(2, \mathcal O)$ of index
     $gq$.  In particular, if the rank of $S/[S,S]$ is one or more, then
   $S$ contains  a  non-congruence subgroup of $PSL(2,\mathcal O)$ of finite index.}

   \smallskip
   \noindent {\em Proof.} Set $G=PSL(2,\mathcal O)$.
Consider the prime decomposition $(q) =\prod \pi_i $,
     $\pi_i \ne \pi_j$ if and only if $i \ne j$.
      Since $\mathcal O\ne
    \Bbb Z$, $ \mathcal O/\pi_i$ is a field of at least 4 elements.
    It follows that $PSL(2, \mathcal O/\pi_i )$ is non-abelian simple,
       $$gcd\,( q,  |SL(2, \mathcal O/n)|/g)= 1  \eqno (5.1)$$
      and that $S$ has a normal subgroup  $N$ of index $q$. It is clear that the
       level of $N$ is a divisor of
       $nq$.
       Suppose that $N$ is congruence.
         By Theorem 3.6, $G(nq) \subseteq N$
         and $G(n) \subseteq S$, where $G(m)$ is the principal congruence subgroup of level
          $m$. We have two cases to consider.

\noindent Case 1.  $G(n)\subseteq N$. It follows that
$N/G(n) $ is a  subgroup of $G/G(n)\cong PSL(2, \mathcal O/n)$ of index $q$. An easy study of the indices of the groups
$G(n) \subseteq N \subseteq S \subseteq G$ implies that $q$ is a
divisor of $|SL(2, \mathcal O /n)|/g$. A contradiction.
Hence $N$ is non-congruence.

\noindent Case 2. $G(n)$ is not a subgroup of $N$. This implies
that $(N\cap G(n))/G(nq)$ is a normal subgroup of
 $ G(n)/G(nq)
\cong SL(2, \mathcal O/q)$
of index $q$.
         We consider the composition factors  of the following normal
      series.
     $$1\triangleleft (N\cap G(n))/G(nq) \triangleleft G(n)/G(nq)
\cong SL(2,\mathcal O/q).\eqno(5.2)$$
    $PSL(2, \mathcal O/\pi_i)$ is a composition factor of (5.2).
    Since $[G(n)  :N\cap G(n)] = q$ is  a prime,
       $PSL(2, \mathcal O/\pi_i)$ is not a composition factor of  $( N\cap G(n))/G(nq) \triangleleft G(n)/G(nq)$.
       Hence  $PSL(2, \mathcal O/\pi_i)$
        is  a composition factor of  $ 1\triangleleft ( N
\cap G(n))/G(nq) $ for all $i$.
          By (5.2), the order of  $G(n)/G(nq)$ is a multiple of $q  |PSL(2, \mathcal O/q)|$. A contradiction.
Hence $N$ in non-congruence.

               Note that $\cap\, xNx^{-1}$ is normal  in $G$.\qed

     \smallskip
     \noindent {\bf Proposition 5.3.} {\em Let $q\in \Bbb N$
be a prime.
Suppose that $q$ is inert or $q\ge 5$ is split in $O_d$  and that $d \ne -1,-3$.
       Then $B_d$
       has a normal non-congruence subgroup of level $q$, index $q$.
       Similarly,  $PE_d$
       has a normal non-congruence subgroup of level $q$, index $q$.
       }
\smallskip

\noindent
{\em Proof.}   Let $S= B_d$. Since $d \ne -1, -3$,
 $S/[S,S] = B_d/B_d'$ has positive rank (Theorem 4.1 and
 Table 1).
 Let $ q $ be given as in the proposition let
    $N$ be a normal subgroup of index $q$ in $S= B_d$.
     By Lemma 5.2, $N$ is normal non-congruence in $S=B_d$.
         Note that the above argument works for $PE_d$.\qed

   \smallskip
     \noindent {\bf Proposition 5.4.} {\em  Let  $d=-1,-3$.
       Then $B_d=PE_d $
       has  normal non-congruence subgroups.
            }

\noindent
{\em Proof.} (i) $d=-1$. Let $S= B_d'''$  be  the third commutator subgroup of $B_d$. Then
the level of $S$ is a divisor of  12, $[B_d : S] = 768$,
 $S$ has 384 generators and $S/[S,S]$ is infinite (Theorem 5.3.1 of Fine [F]).
   By Lemma 5.2,  $B_d$ has a normal  non-congruence subgroup.

\noindent (ii) $d=-3$.  Set $H_2 = B_d (\sqrt {-3})$. Applying results of Fine (Theorem 4.4.4 of Fine [F]), $H_2 = B_d''$ and $S = \left < x^2\,:\, x\in H_2\right >$ is characteristic of index 3 in $H_2$. Further,
$$S =
 \left < w_1, w_2, w_3, w_4
 \,:\, (w_4w_2)^3 = (w_1w_3)^3= (w_1w_2 w_3w_4)^3 = [w_1, w_2]=[w_3,w_4]=1\right >.$$
 Note  that $S/[S,S]\cong \Bbb Z\times \Bbb Z\times \Bbb Z_3\times \Bbb Z_3$.
  By Lemma 5.2,
 $B_d$
       has a  normal  non-congruence subgroup.\qed




\section{Second Application : The power subgroups}
\subsection{}Now we turn our attention to the
  power group
  $B_d^2 = \left < x^2 \,:\, x\in B_d \right > $.
 We shall first prove a lemma that will be useful for our study.
  Throughout the section, $d<0$ is square free,
   $B_d(x)$ is the principal congruence subgroup associated to $x
   $. The main purpose of this section is to show that the number of $d$'s
    such that $B_d^2$ is congruence is finite (Propositions 6.2 and 6.3).

\smallskip
\noindent {\bf Lemma 6.1.} {\em
Let $G = PSL(2, O_d/2)$ and $G^2 = \left < x^2\,:\, x \in G\right >$. Suppose that $(2)$ is a prime ideal of
$O_d$. Then $G=G^2$.
Suppose that $(2)$ is not a prime ideal of $O_d$.
 Then $G/G^2 \cong \Bbb Z_2 \times \Bbb Z_2$.}

\smallskip
\noindent {\em Proof.}
Suppose that $(2)$ is a prime ideal. Then $PSL(2,O/2) \cong A_5$
is non-abelian simple. Hence $G=G^2$.
 In the case $(2)$ is split, $G \cong S_3 \times S_3$ and
 the lemma holds. We may therefore assume that $(2)$ is ramified.
 Let $\pi$ be the
  prime ideal  that divides $(2)=\pi^2$. It follows that $|G|=48$ and
$G$
is an extension  of $B_d(\pi)/B_d(2)=A
\cong \Bbb Z_2 \times \Bbb Z_2 \times \Bbb Z_2$ by $ B_d/B_d(\pi) = X\cong  S_3$.
 To be more accurate,
 pick $w \in \pi \setminus (2)$, then
 $A=\left < x,y,z\right > $, where
  $$
      x=  \left (
\begin{array}{rc}
1 & w \\
0 & 1 \\
\end{array}
\right ),
 y=  \left (
\begin{array}{rc}
1 & 0 \\
w & 1 \\
\end{array}
\right ),
z= \left (
\begin{array}{rc}
1+w & w \\
w & 1+w \\
\end{array}
\right )  .
  \eqno(6.1)
$$
One can show easily
 that the group generated by $x,y$ and $z$ is elementary abelian of order 8
  modulo $(2)$
   (see (i) and (ii) of Lemma 2.4). Similarly, $X= \left < S, R \,\right > , $ where
 $$
  S=  \left (
\begin{array}{rc}
0 & 1 \\
-1 & 0 \\
\end{array}
\right ),
 R=  \left (
\begin{array}{rc}
0 & 1 \\
-1 & -1 \\
\end{array}
\right ).\eqno(6.2)
$$
Note that $S$ and $R$ generate a symmetric group of order 6 modulo $(2)$.
  Let $B$ be the subgroup of $A$ generated by $xz$ and $yz$.
  $$ B = \left <
  xz =
    \left (
\begin{array}{rc}
1+w & 0 \\
w & 1+w \\
\end{array}
\right ),
 yz= \left (
\begin{array}{rc}
1 +w& w \\
0 & 1+w \\
\end{array}
\right ) \right >.\eqno(6.3)$$
$B$ is invariant under the conjugation of $R$.
Hence $B$ and $ R$ generate a subgroup of $G$ of
 order 12.
It is clear that $ \left <z, S\right >\cap\left < B, R\, \right >=1$ and that
 $z$ and $S$ normalise $ D= \left < B, R\,\right >$. Hence $D$ is a normal subgroup of $G$
 and that
 $$G/D \cong  \left <z, S \,\right >\cong \Bbb Z_2 \times \Bbb Z_2.\eqno(6.4)$$
  This implies that $G^2 \subseteq D$. Suppose that $G^2\ne D$. Since $G/G^2$ is an elementary
   abelian 2-group and $|G|=48$, the order of $G^2$ is either 3 or 6. In both cases, one concludes that
   the Sylow 3-subgroup of $G^2$ is a normal subgroup of $G$. This is a contradiction as
   $\left <R\right >$ is not normal.
Hence $D= G^2$. This completes the proof of the lemma.\qed

 \smallskip
\noindent {\bf Proposition 6.2.} {\em
Suppose that $ d \equiv 5\,\,(mod\, 8)$.
Then $B_d'$ and $
  B_d ^2$ are non-congruence if and only if $d \ne -3$.}

 \smallskip
 \noindent {\em Proof.}
 Since  $ d \equiv 5\,\,(mod\, 8)$,
  $(2)\subseteq O_d$ is a prime ideal. It follows that  $PSL(2, O_d/2)$
   $ \cong A_5$.

 \noindent
  (i) $ d\ne -3$.
     $B_d ^2$ is of level 2.
              Suppose that $B_d ^2$ is congruence. By Theorem 3.6,
               $B_d(2) \subseteq B_d ^2$.
Hence
           $  B_d^2 /B_d(2)$ is  a normal subgroup of $B_d/B_d(2)\cong PSL(2, O_d/2)$.
By Lemma 6.1,      $B_d^2 = B_d$.
              This is a contradiction (Theorem 4.1 and Table 1).
           Hence $B_d'$ and $
  B_d ^2$ are not
   congruence.
   \smallskip

   \noindent (ii) $d = -3$.
    Let $w = \sqrt {-3}$.
     $B_{-3}/B_{-3}( w) \cong A_4$  has a unique normal subgroup  $U/B_{-3}(w)$ of
      index 3. In particular, $B_{-3}/U$ is abelain. Hence $B_{-3}'\subseteq U$.
       By Table 1,  $ B_{-3}/B_{-3}' \cong \Bbb Z_3$.
        This implies that $U$ and $ B_{-3}'$ have the same index.
        Hence $ U = B_{-3}'$.
      As a consequence,   $B_{-3}(w) \subseteq U = B_{-3}'$.
       Hence  $ B_{-3}' $ is congruence of level 3. Note that $B_{-3}^2 = B_{-3}$.\qed

\smallskip
\noindent {\bf Proposition 6.3.} {\em
Suppose that
$ d \not \equiv 5\,\,(mod\, 8)$. Then
 $ B_d ^2$ is non-congruence if and only if $|B_d/B_d^2| \ge 8$.
  In particular, if the class number of $O_d$ is three or more, then
   $B_d^2$ is non-congruence.
  }

\smallskip
\noindent {\em Proof.} The level of $B_d^2$ is $2$.
Suppose that $B_d^2$ is congruence. By Theorem 3.6
, $B_d(2) \subseteq B_d^2$.
 Note that $B_d/B_d(2) \cong PSL(2,O_d/2)$.
Since $d\not \equiv 5$ (mod 5), $(2)$ is not a prime.
  By Lemma 6.1,  $B_d/B_d^2 \cong
\Bbb Z_2 \times \Bbb Z_2 $. In particular, $|B_d/B_d^2|  <8$.
 Conversely, suppose that  $|B_d/B_d^2|  <8$. Then
 $|B_d/B_d^2| = 1,2$, or 4. Let $U\subseteq PSL(2, O_d)=B_d $ be the pre-image of $D$
  (see (6.4) of  Lemma 6.1). Then $B_d/U \cong \Bbb Z_2 \times \Bbb Z_2$.
   This implies that $B_d^2 = U$. Since $B_d(2) \subseteq U$, $B_d^2$ is
    congruence.
This completes the proof of the first part of the proposition.

  Suppose that the class number of $O_d$ is three or more.
   By Theorem 4.1, $|B_d/B_d^2| \ge 8$.
    Suppose that $B_d^2$ is congruence. Then $B_d(2) \subseteq B_d^2$ (Theorem 3.6)
     and $B_d^2/B_d(2) \triangleleft  B_d/B_d(2) \cong PSL(2, O_d/2)$.
     In the case $d\equiv 5$ (mod 8), $PSL(2, O_d/2)\cong A_5$ is simple.
     Hence $B_d^2 = B_d$ or $B_d(2)$. This contradicts the fact that
     $B_d/B_d^2$ is abelian of order 8 or more.
      Hence $B_d^2$ is non-congruence.
       In the case  $d\not \equiv 5$ (mod 8), $B_d^2$ is not congruence by the first
        part of the present proposition.\qed

\smallskip
 \noindent {\bf Discussion  6.4.}
  The set of $d$'s $(d<0$ is square free) such that  the class number of $O_d$  is 2 or less is $\{-1,-2,-3,-7,-11,-19,-43,-67,-163\} \cup
  \{ -5,-6,-10,-13,-15,-22,-35,-37,$
  $-51,-58,-91,-115,$
  $
  -123,-187,-235,-267,
  -403,-427\}$.


\subsection{}
  One may apply the proof of Proposition 6.2
  and conclude
immediately that

 \smallskip
\noindent {\bf Proposition 6.5.} {\em
Suppose that $ d \equiv 5\,\,(mod\, 8)$.
Then $PE_d'$ and $
  PE_d ^2$ are non-congruence if and only if $d \ne -3$.}
\smallskip

In the case $d\not\equiv 5$ (mod 8), one has the following.

\smallskip
\noindent {\bf Proposition 6.6.} {\em
Suppose that $ d \not \equiv 5\,\,(mod\, 8)$.
Then  $
  PE_d ^2$ is congruence.}

\smallskip
\noindent {\em Proof.} By Table 1, $PE_d/PE_d^2 \cong \Bbb Z_2 \times \Bbb Z_2$.
 By Lemma 6.1, $PE_d/PE_d(2) \cong PSL(2, O_d/2)$ has a normal subgroup $D$
  such that $PSL(2, O_d/2)/D \cong \Bbb Z_2 \times \Bbb Z_2$. Denoted by $V$
   the pre-image of $D$ in $PE_d$. Then $PE_d/V \cong \Bbb Z_2 \times \Bbb Z_2$.
    Hence $PE_d^2 = V$.
    Note that $PE_d(2) \subseteq V$.
    Hence $PE_d^2$  is congruence. \qed

\section{Third Application : Congruence Subgroup Property (CSP)}
The main purpose of this section is to show that if $G = PSL(2, O_d)$ has CSP, then
 $G/[G,G]$ is a finite group of order $2^a 3^b$ for some $a, b$ (Proposition 7.2).
Note that $B_d$ do not have CSP and $B_d/[B_d, B_d]$ is finite only if
 $d = -1,  -3$.

 \smallskip
 \noindent {\bf Proposition 7.1.} {\em  Suppose that $G =PSL(2, \mathcal O)  $
  has CSP.
   Let $S$ be a  subgroup of $G  $ of index $g$, level $n$.
     Then $|S/[S, S]|$ is finite. Let $q$ be a prime divisor of
   $|S/[S, S]|$. Suppose that
     $q$ is inert or $q\ge 5$ is split. Then
     $q$ divides
      $ |SL(2, \mathcal O/n)|/g$.
    In particular, if $q$ is a prime divisor of   $|G/[G,G]| <\infty$,  then $q$ is either ramified or split. In the case $q$ is split, $q=2$ or $3$.
 }

    \smallskip
     \noindent {\em Proof.} Since $G$ has CSP, $S/[S,S]$ is finite (Section 16 of [BMS]). The second part of the proposition
follows immediately by applying
Lemma 5.2. \qed

 \smallskip

     \noindent {\bf Proposition  7.2.}
     {\em
      Let $O_d$ be the ring of integers of the real quadratic  field $\Bbb Q(\sqrt d\,)$
      and let  $G =PSL(2, O_d)$. Let $q\in \Bbb N$ be a prime.
       Suppose that $q$ divides $|G/[G, G]|$. Then $q$ is ramified or split
      and  $q = 2$ or $3$.  In particular,
      $|G/[G, G]| = 2^a 3^b$.
      }

\smallskip
\noindent {\em Proof.} Since $\Bbb Q[\sqrt d\,]$ has a unit
of infinite order, $G$ has CSP (Serre [S]).
By our assumption, $G$ has a normal subgroup of index $q$.
 By Lemma A1, $q$ is ramified or split and $q\le 5$.\qed

\smallskip
Proposition 7.2 suggests the following.

 \smallskip

     \noindent {\bf Conjecture 7.3.}
     {\em  Suppose that $G = PSL(2, \mathcal O)$ has CSP.
       Let  $q$ be a prime divisor of  $|G/[G, G]|$. Then $q$ is ramified or split
      and  $q = 2$ or $3$.  In particular,
      $|G/[G, G]| = 2^a 3^b$.
      }


\section {Appendix A}

Let $q\ge 5$ be a rational prime that is ramified  in $O_d$, the ring of integers of the quadratic field $\Bbb Q[\sqrt d]$, $d \in \Bbb Z$.
  Let $\pi$ be the prime ideal that
divides $(q)=\pi^2$ and let  $\{q, x\}$ be an integral basis of $\pi$.
Let $G= PSL(2, O_d)$. Then $G/G(q)$
is an extension of $G(\pi)/G(q)$ by $G/G(\pi)
 \cong PSL(2, O_d/\pi)=PSL(2,q)$,
where $G(\pi)/G(q)
 \cong \Bbb Z_q \times
 \Bbb Z_q \times
 \Bbb Z_q .$
  $$
G(\pi)/G(q)\cong  \left <
  r =  \left (
\begin{array}{rc}
1 & 0 \\
x & 1 \\
\end{array}
\right ),
 s=  \left (
\begin{array}{rc}
1 & x \\
0 & 1 \\
\end{array}
\right ),
u= \left (
\begin{array}{rc}
1-x & x \\
-x & 1+x \\
\end{array}
\right ) \right > .
  \eqno(A1)
$$
Let $t= rs^{-1}u$.
 $t$ (mod $q$) can be found in $(A2)$. It is clear that $\left < r,s,u\right > = \left <r,s,t\right >$ (mod $q$).
 $$
 G/G(\pi) \cong
  \left <
  S =  \left (
\begin{array}{rc}
0 & 1 \\
-1 & 0 \\
\end{array}
\right ),
 T=  \left (
\begin{array}{rc}
1 & 1 \\
0 & 1 \\
\end{array}
\right )
 \right >,\,
t\equiv \left (
\begin{array}{rc}
1-x & 0 \\
0 & 1+x \\
\end{array}
\right ).
  \eqno(A2)$$
Set $v^g = gvg^{-1}$. Direct calculation shows that (see (i) and (ii) of Lemma 2.4 for
matrix multiplication of $G(\pi)/G(q)$)
$$r^S=s^{-1}, s^S= r^{-1}, t^S= t^{-1},
r^T = rs^{-1} t^{-1}, s^T = s, t^T = s^2t.\eqno(A3)$$

\noindent {\bf Lemma A1.} {\em
Suppose that $G = PSL(2, O_d)$ has CSP.
Let $q$ be a rational
prime. Suppose that $q\in O_d$ is inert or $q\ge 5$.
 Then
 $G $ has no normal subgroups of index $q$.}

\smallskip
\noindent {\em Proof.}
Suppose that $G$ has a normal subgroup $K$ of index $q$. Since $G$ has CSP, $G(q)\subseteq K$ (Theorem 3.6).
In the case $q$ is inert, $G/K$ is a normal subgroup of the
simple group $G/G(q)\cong PSL(2, q^2)$. A contradiction. In the case $q\ge 5$ is
split, $G/K$ is a normal subgroup of index $q$ of $G/G(q) \cong (SL(2, q)\times SL(2,q))
/\Bbb Z_2$. Again, a contradiction. Hence we may assume that $q\ge 5$ is ramified.
Let $\pi$ be the  prime ideal that divides $(q)=\pi^2$. Since $q\ge 5$,
 $PSL(2, O_d/\pi)$ is non-abelian simple. Hence $G(\pi)$ is not
a subgroup of $K$.
Hence $N = K/G(q)  \cap G(\pi)/G(q)$ is a normal subgroup of $G/G(q)$ of
order $q^2$.
Take $1 \ne\sigma = r^as^bt^c \in N\subseteq G(\pi)/G(p)$. There are four cases to consider.

\noindent {\bf Case 1.} $a \ne b$.  Recall that $N$ is normal.
$
\sigma^S\sigma
=r^{a-b}s^{b-a} \in N
$.
 Hence $rs^{-1} \in N$. It follows that
$(rs^{-1})^{-1}(rs^{-1})^T =( st)^{-1} \in N$. Hence
$(st)(st)^{-1})^T =s^{-2} \in N$. It follows that $s\in N$. By (A3) and the fact
 that $N$ is normal, $r, t \in N$. Hence $N$ has order $q^3$. A contradiction.

\noindent {\bf Case 2.} $a=b=0$. Hence $t^c \in N$. This implies
that $t\in N$. By (A3) and the fact that $N$ is normal, $r,s \in N$. Hence $N$ has order $q^3$.
A contradiction.

\noindent {\bf Case 3.} $a=b \ne 0$ and $2c/a-1 =0$ (mod $q$).
Hence $r^as^at^c\in N$. It follows that $rst^d \in N$, where
 $d= c/a$. Hence $(rst^d)^{-1}(rst^d)^T=t^{-1} \in N$. Similar
to Case 2, this is a contradiction.

\noindent {\bf Case 4.} $a=b\ne 0$ and $2c/a -1 \ne 0$ (mod $q$).
Similar to Case 3, $rst^d\in N$, where $d=c/a$.
Hence $(rst^d)^{-1}(rst^d)^T =s^{2d-1}t^{-1} \in N$. It follows that $(s^{2d-1}t^{-1})^S(s^{2d-1}t^{-1})=r^{1-2d}s^{2d-1}\in N$.
Hence $rs^{-1}\in N$. Similar to Case 1, this is a contradiction.

\smallskip
Hence $G$ has no normal subgroups of index $q$. \qed


\medskip
{\small
\noindent Cheng Lien Lang\\
\noindent Department of Mathematics, I-Shou  University, Kaohsiung, Taiwan,
Republic of China.

\noindent   \texttt{cllang@isu.edu.tw}

\medskip
\noindent Mong Lung Lang\\
\noindent Singapore 669608, Republic of Singapore.

\noindent \texttt{lang2to46@gmail.com}}

\medskip


\end{document}